\documentclass[12pt]{article}
 \usepackage{amsmath,amssymb,amsthm}
 \RequirePackage[dvips]{graphicx}
 \def\draw #1 by #2 (#3){
  \vbox to #2{
    \hrule width #1 height 0pt depth 0pt
    \vfill
    \special{picture #3} 
    }
  }

 \def\scaleddraw #1 by #2 (#3 scaled #4){{
  \dimen0=#1 \dimen1=#2
  \divide\dimen0 by 1000 \multiply\dimen0 by #4
  \divide\dimen1 by 1000 \multiply\dimen1 by #4
  \draw \dimen0 by \dimen1 (#3 scaled #4)}
  }

\newtheorem{theorem}{Theorem}[section]
\newtheorem{example}{Example}
\newtheorem{problem}[example]{Problem}
\newtheorem{defin}[theorem]{Definition}
\newtheorem{lemma}[theorem]{Lemma}

\newtheorem{remark}[theorem]{Remark}

\newtheorem{nt}{Note}

\newenvironment{pf}{\medskip\noindent{Proof:  \hspace*{-.4cm}}
       \enspace}{\hfill \qed \newline \medskip}

 \newcommand{\singlespacing}{\let\CS=\@currsize\renewcommand{\baselinestretch}{1}\tiny\CS}
 \newcommand{\oneandahalfspacing}{\let\CS=\@currsize\renewcommand{\baselinestretch}{1.25}\tiny\CS}
 \newcommand{\doublespacing}{\let\CS=\@currsize\renewcommand{\baselinestretch}{1.35}\tiny\CS}

 \newtheorem{rule-def}[theorem]{Rule}

\begin{document}
\baselineskip 16pt
 \newcommand{\la}{\lambda}
 \newcommand{\si}{\sigma}
 \newcommand{\ol}{1-\lambda}
 \newcommand{\be}{\begin{equation}}
 \newcommand{\ee}{\end{equation}}
 \newcommand{\bea}{\begin{eqnarray}}
 \newcommand{\eea}{\end{eqnarray}}

 \begin{center}
 {\Large \bf Open problem on $\sigma$-invariant}\\

  \vspace{10mm}

 {\large \bf Kinkar Ch. Das$^a$} and {\large \bf Seyed Ahmad Mojallal$^b$}

 \vspace{9mm}

 \baselineskip=0.20in

 {\footnotesize $^a${\it Department of Mathematics, Sungkyunkwan University, Suwon 16419, Rep. of Korea}}\\
 {\footnotesize $^b${\it Applied Algebra and Optimization Research Center, Sungkyunkwan University,}\\ Suwon 16419, Rep. of Korea,\/} \\
  {\tt kinkardas2003@googlemail.com}, {\tt mojallal@skku.edu}

 \vspace{4mm}

 \end{center}

 \vspace{5mm}

 \baselineskip=0.23in

 \begin{abstract} Let $G$ be a graph of order $n$ with $m$ edges. Also let $\mu_1\geq \mu_2\geq \cdots\geq \mu_{n-1}\geq \mu_n=0$ be the Laplacian eigenvalues
 of graph $G$ and let $\sigma=\sigma(G)$ $(1\leq \sigma\leq n)$ be the largest positive integer such that $\mu_{\sigma}\geq \frac{2m}{n}$.
 In this paper, we prove that $\mu_2(G)\geq \frac{2m}{n}$ for almost all graphs. Moreover, we characterize the extremal graphs for any graphs.
 Finally, we provide the answer to Problem 3 in \cite{KMT}, that is, the characterization of all graphs with $\sigma=1$.

 \bigskip

 \noindent
 {\bf Key Words:} Graph, Laplacian matrix, Second largest Laplacian eigenvalue, Average degree, Laplacian energy, $\sigma$-invariant\\
 \\
 {\bf 2000 Mathematics Subject Classification:} 05C50
 \end{abstract}

 \baselineskip=0.30in

 \section{Introduction}

 Let $G =(V,\,E)$ be a simple graph with vertex set $V(G)=\{v_1,\,v_2,\ldots,\,v_n\}$ and edge set $E(G)$, where $|V(G)|=n$, $|E(G)|=m$. Let
 $d_i$ be the degree of the vertex $v_i$ $(i=1,\,2,\ldots,\,n)$. The maximum vertex degree is denoted by $\Delta_1\,(=d_1)$
 and the second maximum by $\Delta_2\,(=d_2)$. Let $N(v_i)$ is the neighbor set of vertex $v_i$, $i=1,\,2,\ldots,\,n$. Let $A(G)$ and $D(G)$ be the  adjacency matrix and the diagonal
 matrix of vertex degrees of $G$, respectively. The Laplacian matrix of $G$ is  $L(G)=D(G)-A(G)$. This matrix has nonnegative eigenvalues $n \geq \mu_1 \geq \mu_2 \geq \cdots \geq \mu_n=0$.
 When more than one graph is under consideration, then we write $\mu_i(G)$ instead of $\mu_i$.

 \vspace*{3mm}

 For a graph $G$, consider the positive number $\sigma=\sigma(G)$ $(1\leq \sigma\leq n)$ of the Laplacian eigenvalues greater than or equal to the average degree $\frac{2m}{n}$.
 More precisely $\sigma$ is the largest positive integer for which $\mu_{\sigma}\ge \frac{2m}{n}$. This number as a spectral graph invariant with several open problems and conjectures are
 introduced in \cite{KMT}.

 \vspace{3mm}

 Let $I$ be an interval of the real line. Denote by $m_G(I)$ the number of Laplacian eigenvalues, multiplicities included, that belong to $I$.
 Notice that $m_G(I)$ is a natural extension of multiplicity $m_G(\mu)$ of a Laplacian eigenvalue $\mu$.  Merris in \cite{MERIS} presented several results on $m_G(I)$ and gave some references
 for its applications. This research topic was extensively investigated in many papers (see, \cite{GMS,MW,HJT,Me2,MERIS}). Indeed, by the definition of $\sigma$,
 we have $$\sigma(G)=m_G\Big(\Big[\frac{2m}{n},\,n\Big]\Big).$$
 It is worth noticing that the value of $\sigma(G)$ sheds light on the distribution of the Laplacian eigenvalues of a graph $G$. Actually it determines how many
 Laplacian eigenvalues of graph $G$ are greater than or equal to the average of the Laplacian eigenvalues of graph.

 \vspace{3mm}

  A further Laplacian--spectrum--based graph invariant was put forward by Gutman and Zhou \cite{GuZh} as
 \begin{equation}                    \nonumber
 LE=LE(G)=\sum^n_{i=1}\left|\mu_i-\frac{2m}{n}\right|.
 \end{equation}
 For its basic properties, including various lower and upper bounds, see \cite{DA2,KI-MO,L15}. It is not difficult to see that
 \begin{equation}
 LE(G)=2\sum^{\sigma}_{i=1}\mu_i-\frac{4m \sigma}{n}.\label{ghs1}
 \end{equation}
 This is another motivation to study this graph invariant for the Laplacian energy of a graph $G$.

 \vspace{3mm}

 Therefore, the spectral parameter $\sigma$ is reasonable relevant in spectral graph theory. To know more information about this spectral graph invariant
 and its applications, see \cite{KMG, KMT, PG}. In \cite{KMG}, all graphs with $\sigma(G)=n-1$ were characterized and the result was applied for Laplacian
 energy of graphs. It is interesting problem to characterize all graphs for some specific value of $\sigma=\sigma(G)$ between 1 and
 $n-2$. In particular, the following problem is given in \cite{KMT}:
 \begin{problem}\label{prob3}{\rm \cite{KMT}} Characterize the graphs with $\sigma=1$.
 \end{problem}

 \vspace{3mm}

 Li and Pan \cite{LP} showed that
 \begin{equation}
 \mu_2(G)\geq \Delta_2\label{1ekine1}
 \end{equation}
 with equality if $G$ is an $r\times s$ complete bipartite graph $K_{r,\,s}$ $(r+s=n)$ or a tree $T$ with degree sequence
 $\pi(T)=(n/2,\,n/2,\underbrace{1,\ldots,\,1}_{n-2})$, where $n\geq 4$ is even. This result was improved by one of the present authors in \cite{DAS} as follows:
 \begin{equation} \nonumber
 \mu_2(G) \geq \left\{
    \begin{array}{ll}
      \frac{\Delta_2+2+\sqrt{(\Delta_2-2)^2+4c_{12}}}{2} & ~~\hbox{if } v_1v_2\in E(G) \\
          &   \\
      \frac{\Delta_2+1+\sqrt{(\Delta_2+1)^2-4c_{12}}}{2} & ~~\hbox{if } v_1v_2\notin E(G)\,,
    \end{array}
  \right.
 \end{equation}
 where $v_1$ and $v_2$ are the maximum and the second maximum degree vertices of graph $G$, respectively, and $c_{12}=|N(v_1)\cap N(v_2)|$.
 We refer to  \cite{DAS,LP,WYS} for more background on the second largest Laplacian eigenvalues of graphs.

 \vspace*{3mm}

 For two vertex-disjoint graphs $G_1$ and $G_2$, We use $G_1 \cup G_2$ to denote their union.
 The join $G_1\vee G_2$ of graphs $G_1$ and $G_2$ is the graph obtained from the disjoint union of $G_1$ and $G_2$ by adding all edges between $V(G_1)$ and $V(G_2)$.
 For any two sets $A,\,B\subseteq V(G)\,(E(G))$, let $A\backslash B$ be the set of vertices (edges) belongs to $A$, but not
 $B$. Denoted by $|A|$, is the cardinality of the set $A$. As usual, $K_n$, $K_{1,n-1}$ and $DS_{p,\,q}$ $(p\geq q\geq 2,\,p+q=n)$, denote, respectively, the complete
 graph, the star graph, and the double star graph on $n$ vertices.

 \vspace*{3mm}

 The paper is organized as follows. In Section 2, we give a list of some previously known results. In Section 3, we give a lower bound on the second largest Laplacian eigenvalue of graph $G$
 and characterize the extremal graphs. We present an upper bound for the third smallest Laplacian eigenvalue of $G$. In Section 4, we obtain the solution for
 Problem \ref{prob3}. Finally we apply this result for Laplacian energy of graphs.

 \section{Preliminaries}

 In this section, we shall list some previously known results that will be needed in the next two sections.
 We begin with first two results on symmetric matrices of order $n$.
 \begin{lemma} \label{t1} {\rm \cite{FA}} Let $A$ and $B$ be two real symmetric matrices of size $n$. Then for any $1\leq k\leq n$,
   $$\sum^k_{i=1}\lambda_i(A+B)\leq \sum^k_{i=1}\lambda_i(A)+\sum^k_{i=1}\lambda_i(B),$$
  where $\lambda_i(M)$ is the $i$-th largest eigenvalue of $M$ $(M=A,\,B)$.
 \end{lemma}

 \begin{lemma} {\rm \cite{SC}} \label{k0} Let $B$ be an $n\times n$ symmetric matrix and let $B_k$ be its leading $k\times k$ submatrix. Then, for
 $i=1,\,2,\ldots,\,k$,
 \begin{equation}                    \label{kh1}
 \lambda_{n-i+1}(B) \leq \lambda_{k-i+1}(B_k) \leq \lambda_{k-i+1}(B),
 \end{equation}
 where $\lambda_i(B)$ is the $i$-th largest eigenvalue of $B$.
 \end{lemma}

 \noindent
 The following result is well-known as interlacing theorem on Laplacian eigenvalues.
 \begin{lemma} {\rm\cite{HE}} \label{m1} Let $G$ be a graph of $n$ vertices and let $H$ be a
 subgraph of $G$ obtained by deleting an edge in $G$. Then
 \begin{equation}
 \mu_1(G)\geq \mu_1(H)\geq \mu_2(G)\geq \mu_2(H)\geq \mu_3(G)\geq\cdots\geq \mu_{n-1}(G)\geq \mu_{n-1}(H)
 \geq \mu_n(G)\geq \mu_n(H)\geq 0,\nonumber
 \end{equation}
 where $\mu_i(F)$ is the $i$-th largest Laplacian eigenvalue of the graph $F$.
 \end{lemma}

 \noindent
 Merris in \cite{MERIS} gave a lower bound on Laplacian spectral radius of a graph $G$ as follows:
 \begin{lemma} {\rm\cite{MERIS}} \label{k6} Let $G$ be a graph on $n$ vertices which has at least one edge. Then
 \begin{equation}
 \mu_1\geq \Delta_1+1.\label{lu1}
 \end{equation}
 Moreover, if $G$ is connected, then the equality holds in (\ref{lu1}) if and only if $\Delta_1=n-1$.
 \end{lemma}

 \noindent
 We now mention an upper bound on the Laplacian spectral radius of graph $G$:
 \begin{lemma} {\rm \cite{ANDER}} \label{kk}  Let $G$ be a graph. Then
  $$\mu_1(G) \leq \max\Big\{d_i+d_j\,|\,(v_i,\,v_j) \in E(G)\Big\}\leq \Delta_1+\Delta_2,$$
 where $d_i$ is the degree of vertex $v_i\in V(G)$.
 \end{lemma}

 \noindent
 The following result is obtained in \cite{DAS3}.
 \begin{lemma} \label{t4} {\rm \cite{DAS3}} Let $G$ be a connected graph with $n\geq 3$ vertices. Then $\mu_2=\mu_3=\cdots=\mu_{n-1}$ if and only if
 $G \cong K_n,\,G\cong K_{1,\,n-1}$ or $G \cong K_{\frac{n}{2},\,\frac{n}{2}}$ $(n\mbox{ is even})$.
 \end{lemma}

 \noindent
 Pan and Hou \cite{PH} obtained the necessary condition for a graph to have an equality of the second largest Laplacian
 eigenvalue $\mu_2$ and its lower bound $\Delta_2$:
 \begin{lemma}  \label{k2} {\rm \cite{PH}} Let $G$ $(\ncong K_{1,\,n-1})$ be a connected graph of order $n\geq 3$ with maximum degree vertex $v_1$ and the second maximum degree vertex $v_2$. If
 $\mu_2(G)= \Delta_2$, then\\
 $(1)$ $N(v_1)\cap N(v_2)=\emptyset$;\\
 $(2)$ $\Delta_1=\Delta_2$;\\
 $(3)$ $\Delta_1+\Delta_2=n$.
 \end{lemma}

 \noindent
 The nice relation between Laplacian spectrum of graph $G$ and the Laplacian spetrum of graph $\overline{G}$ is the following:
 \begin{lemma} {\rm \cite{MERIS}} \label{j2} Let $G$ be a graph with Laplacian spectrum $\{0=\mu_n,\,\mu_{n-1},\ldots,\mu_2,\,\mu_1\}$.
 Then the Laplacian spectrum of $\overline{G}$ is $\{0,\,n-\mu_1,\,n-\mu_2,\ldots,\,n-\mu_{n-2},\,n-\mu_{n-1}$\}, where
 $\overline{G}$ is the complement of the graph $G$.
 \end{lemma}

 \section{Lower bound for the second largest Laplacian eigenvalue of graphs}

 In this section we give some lower bounds on the second largest Laplacian eigenvalue of graph $G$ and characterize the extremal graphs.
 Moreover, we give an upper bound for the third smallest Laplacian eigenvalue of graph
 $G$. From now we always assume that $\Delta_1=d_1\ge \Delta_2=d_2 \ge \cdots \ge d_n$.
 Let $k$ $(2\leq k\leq n)$ be the largest positive integer such that $d_2=\cdots=d_k=\Delta_2$.
 For $G\cong DS_{p,\,q}$ $(p\geq q\geq 2,\,p+q=n)$, we have $\Delta_2=d_2=q>1=d_3$ and hence $k=2$. For $G\cong K_{1,\,n-1}$, we have $d_2=d_3=\cdots=d_n=\Delta_2$ and $k=n$.
 For any graph $G$ we have
 \begin{equation}
 \frac{2m}{n}=\frac{\Delta_1+\Delta_2+\sum\limits^n_{i=3}\,d_i}{n}\leq \Delta_2+\frac{\Delta_1-\Delta_2}{n}<\Delta_2+1.\label{mo4}
 \end{equation}
 Now we have the following result:
 \begin{lemma} \label{1a1} Let $G$ be a graph of order $n>2$ with $m$ edges. If $\Delta_2<\frac{2m}{n}$, then
 \begin{eqnarray}
 (i)&& \Delta_1-\Delta_2> \sum_{i=3}^n (\Delta_2-d_i),  \label{mo1}\\
 (ii)&& \Delta_2=d_3,\nonumber\\
 (iii)&&\Delta_1>\Delta_2+n-k,  \label{mo2}\\
 (iv)&& k\geq \Delta_2+n+1-\Delta_1\geq \Delta_2+2.\label{mo3}
 \end{eqnarray}
 \end{lemma}

 \begin{pf} $(i)$ Since $\Delta_2<\frac{2m}{n}$, we have
 \begin{equation}
 \sum_{i=1}^n \,d_i=2m>n\,\Delta_2,~~\mbox{ that is, }~~\Delta_1-\Delta_2> \sum_{i=3}^n (\Delta_2-d_i). \nonumber
 \end{equation}

 \noindent
 $(ii)$ We have $\Delta_2\geq d_3$. If $\Delta_2>d_3$, then
   $$\Delta_1>\Delta_2+\sum_{i=3}^n (\Delta_2-d_i)\geq \Delta_2+n-2\geq n-1,\mbox{ a contradiction.}$$
 Hence $\Delta_2=d_3$.

 \noindent
 $(iii)$ From $(\ref{mo1})$, we have
   \begin{equation}
   \Delta_1\geq \Delta_2+1+\sum_{i=k+1}^n (\Delta_2-d_i)\geq \Delta_2+1+n-k>\Delta_2+n-k.   \nonumber
   \end{equation}

 \noindent
 $(iv)$ From $(\ref{mo2})$, we have
    \begin{equation}
    k\geq \Delta_2+n+1-\Delta_1\geq \Delta_2+2. \nonumber
    \end{equation}
 This completes the proof of the result.
  \end{pf}

 \noindent
 Let $v_1$ be the maximum degree vertex in $G$, that is, $d_1=\Delta_1$. Let $S$ be the set of vertices $v_j\in V(G)\backslash\{v_1\}$ such that $d_j=\Delta_2$, that is,
        $$S=\Big\{v_j\in V(G)\backslash\{v_1\}\,|\,d_j=\Delta_2\Big\}.$$
 Since $\Delta_2=d_2$, we have $|S|\geq 1$. Denote by $T=V(G)\backslash (S\cup \{v_1\})$ (see, Fig. 1). Then we have $d_i\leq \Delta_2-1$ for any vertex $v_i\in T$.
 Since $d_2=\cdots=d_k=\Delta_2$, we have $|S|=k-1$ and hence $|T|=n-k$. Let $S_1$ and $S_2$ be two sets of vertices such that
     $$S_1=\Big\{v_j\in S\,|\,v_1v_j\in E(G)\Big\},\,S_2=S\backslash S_1.$$

 \begin{center}
 \includegraphics[height=7cm,keepaspectratio]{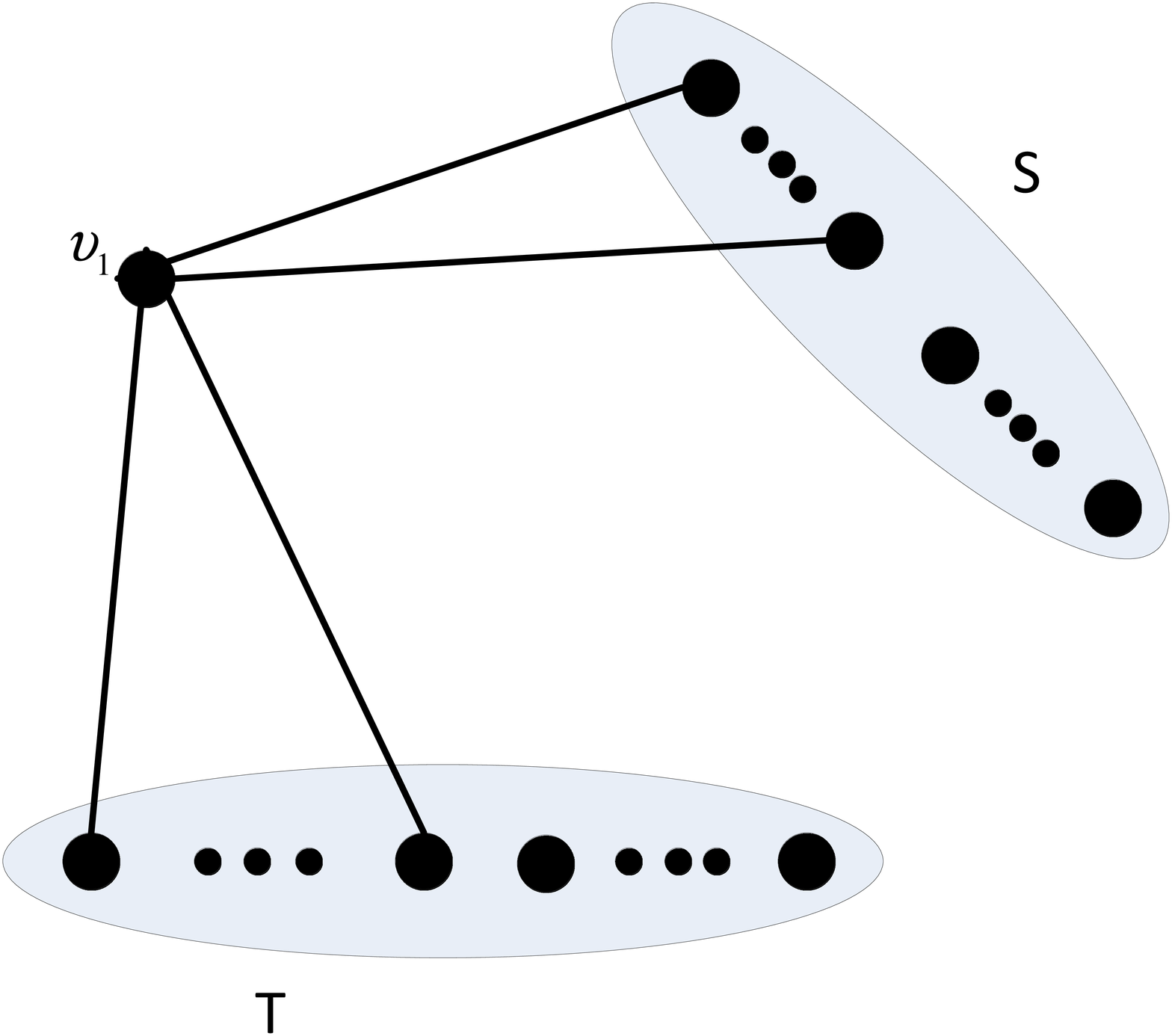}\\
  $G$\\
 \vspace*{3mm}
 \noindent
 Fig. 1. The graph $G$ with the vertex set $V(G)=\{v_1\}\cup S\cup T$\,.
 \end{center}

 We need the following two results for our main result on the lower bound of $\mu_2$.
 \begin{lemma} \label{1a2} Let $G$ be a graph of order $n>2$ with $m$ edges. If there is an edge in the subgraph induced by $S$ in $G$ and
 $\Delta_2<\frac{2m}{n}$, then $\mu_2(G)>\frac{2m}{n}$.
 \end{lemma}

 \begin{pf} Let $v_i$ and $v_j$ be two vertices in the subgraph induced by $S$ in $G$ such that $v_iv_j\in E(G)$.
 Here we consider the following three cases:

 \vspace*{3mm}

 \noindent
 ${\bf Case\,1:}$ $v_i,\,v_j\in S_1$. In this case we consider $3\times 3$ leading principal submatrix $L_1$ corresponding to vertices $v_1$,
 $v_i$ and $v_j$ of
 Laplacian matrix $L(G)$, where
 \[ L_1=\left( \begin{array}{ccc}
 \Delta_1 &-1 &-1  \\
 -1 &\Delta_2 &-1  \\
 -1 &-1 &\Delta_2  \\
 \end{array} \right).\]

 \noindent
 The characteristic polynomial of $L_1$ is
    $$f(x)=x^3-(\Delta_1+2\Delta_2)x^2+(2\Delta_1\,\Delta_2+\Delta_2^2-3)x-\Delta_1\,\Delta^2_2+\Delta_1+2\,\Delta_2+2.$$

 \noindent
 The roots of $f(x)=0$ are
   $$\frac{\Delta_1+\Delta_2-1}{2}\pm \frac{1}{2}\,\sqrt{(\Delta_1-\Delta_2)(\Delta_1-\Delta_2+2)+9}\,,\,\Delta_2+1.$$

 \noindent
 By Lemma \ref{k0} with (\ref{mo4}), we have
    $$\mu_2(G)\geq \mu_2(L_1)=\Delta_2+1>\frac{2m}{n}\,.$$

 \vspace*{3mm}

 \noindent
 ${\bf Case\,2:}$ $v_i\in S_1,\,v_j\in S_2$. Since $\Delta_2<\frac{2m}{n}$, from $(\ref{mo2})$ with $|T|=n-k$, we get $\Delta_1\geq \Delta_2+1+|T|$.
 Since $S\cup T\cup \{v_1\}=V(G)$, therefore there are at least $\Delta_2+1$ vertices in $S$, adjacent to
 $v_1$. From the above, $|S_1|\geq \Delta_2+1\geq 3$ $(\mbox{as }\Delta_2\geq d_i\geq 2)$ and $|S_2|\geq 0$. From this result, we can assume that a vertex set
 $W=\{w_1,\ldots,\,w_{\Delta_2}\} \subseteq S_1$, where $v_i\notin W$. If there exists an edge $(\mbox{say, }v_rv_{\ell})$ in the subgraph induced by $W\cup \{v_i\}$ in $G$, then by {\bf Case 1}
 ($L_1$ corresponding to vertices $v_1$, $v_r$ and $v_{\ell}$), we get the required result.
 Otherwise, $W\cup \{v_i\}$ is an independent set in $G$. First we assume that $v_jw_i\notin E(G)$ for $1\le i \le \Delta_2$. Then $K_{1,\,\Delta_2}\cup K_{1,\,\Delta_2}$ is a subgraph of $G$
 with one star $K_{1,\,\Delta_2}$ induced by $\{v_1,\,w_1,\,w_2,\ldots,\,w_{\Delta_2}\}$ and the other star $K_{1,\,\Delta_2}$ induced by $\{v_j\} \cup
 N(v_j)$, where $\{v_1,\,w_1,\,w_2,\ldots,\,w_{\Delta_2}\}\cap (\{v_j\} \cup N(v_j))=\emptyset$ and $|N(v_j)|=\Delta_2$.
 By Lemma \ref{m1}, we have
     $$\mu_2(G)\geq \mu_2(K_{1,\,\Delta_2}\cup K_{1,\,\Delta_2})=\Delta_2+1>\frac{2m}{n},\,\mbox{ by (\ref{mo4})}.$$

 Next we assume that there exists at least one vertex in $W$, $w_q$ (say), adjacent to $v_j$. In this case we consider $4\times 4$ leading principal submatrix $L_2$ corresponding to vertices
 $v_1$, $v_i$, $w_q$ and $v_j$ of Laplacian matrix $L(G)$, where
 \[ L_2=\left( \begin{array}{cccc}
 \Delta_1 &-1 & -1&0  \\
 -1 &\Delta_2 & 0&-1  \\
 -1&0 &\Delta_2 &-1 \\
 0&-1 &-1 &\Delta_2 \\
 \end{array} \right).\]

 \noindent
 The characteristic polynomial of $L_2$ is
   $$g(x)=x^4 - (\Delta_1+ 3 \Delta_2)x^3+(3\Delta_1 \Delta_2 + 3 \Delta_2^2 - 4)x^2-(3\Delta_1 \Delta_2^2 + \Delta_2^3 - 2\Delta_1-6\Delta_2)x+\Delta_1 \Delta_2^3-2 \Delta_1 \Delta_2-2 \Delta_2^2.$$

 \noindent
 We have
 \begin{equation}
 g(\Delta_2+1)=\Delta_1-\Delta_2-3~~~~\mbox{and}~~~g(\Delta_1)=-2(\Delta_1-\Delta_2)^2<0~\mbox{ as }~\Delta_1>\Delta_2.\label{1km1}
 \end{equation}
 From $\Delta_1\geq \Delta_2+1$, we now consider the following three subcases:

 \vspace*{2mm}

 \noindent
 ${\bf Subcase\,2.1:}$ $\Delta_1=\Delta_2+1$. Thus we have
 $$g\left(\Delta_2+\frac{1}{3}\right)=\frac{17}{27} \Delta_1 - \frac{17}{27} \Delta_2 - \frac{35}{81}=\frac{16}{81}>0.$$
 By Lemma \ref{k6} with the above result and (\ref{1km1}), we conclude that $\mu_2(L_2)> \Delta_2+\frac{1}{3}$. Since in this case $n\ge 4$,
 by Lemma \ref{k0}  and (\ref{mo4}), we get $$\mu_2(G)\ge \mu_2(L_2)> \Delta_2+\frac{1}{3}> \Delta_2+\frac{1}{n}=\Delta_2+\frac{\Delta_1-\Delta_2}{n}\ge \frac{2m}{n}.$$

 \vspace*{2mm}

 \noindent
 ${\bf Subcase\,2.2:}$ $\Delta_1=\Delta_2+2$. We have $$g\left(\Delta_2+\frac{2}{3}\right)=\frac{28}{27} \Delta_1-\frac{28}{27} \Delta_2 - \frac{128}{81}=\frac{40}{81}>0.$$
 Again by Lemma \ref{k6} with the above result and (\ref{1km1}), we have $\mu_2(L_2)> \Delta_2+\frac{2}{3}$. Since in this case $n\ge 4$, again by Lemma \ref{k0}  and (\ref{mo4}), we get
 $$\mu_2(G)\ge \mu_2(L_2)> \Delta_2+\frac{2}{3}> \Delta_2+\frac{2}{n}=\Delta_2+\frac{\Delta_1-\Delta_2}{n}\ge \frac{2m}{n}.$$

 \vspace*{2mm}

 \noindent
 ${\bf Subcase\,2.3:}$ $\Delta_1\ge \Delta_2+3$. We have $g(\Delta_2+1)=\Delta_1-\Delta_2-3\ge 0$. Similarly, as the above subcases, we have $\mu_2(L_2)> \Delta_2+1$.
 Again by Lemma \ref{k0} with (\ref{mo4}), we have $$\mu_2(G)\geq \mu_2(L_2)\ge \Delta_2+1>\frac{2m}{n}\,.$$

 \vspace*{3mm}

 \noindent
 ${\bf Case\,3:}$ $v_i,\,v_j\in S_2$.  In this case we consider $3\times 3$ leading principal submatrix $L_3$
 corresponding to vertices $v_1$, $v_i$ and $v_j$ of
 Laplacian matrix $L(G)$, where
 \[ L_3=\left( \begin{array}{ccc}
 \Delta_1 &0 &0  \\
 0 &\Delta_2 &-1  \\
 0 &-1 &\Delta_2  \\
 \end{array} \right).\]

 \noindent
 The eigenvalues of $L_3$ are $\Delta_1$, $\Delta_2+1$ and $\Delta_2-1$. Since $\Delta_2<\frac{2m}{n}$, (\ref{mo2}) holds, that is, $\Delta_1\geq \Delta_2+1$. Again by Lemma \ref{k0} with (\ref{mo4}), we have
    $$\mu_2(G)\geq \mu_2(L_3)=\Delta_2+1>\frac{2m}{n}.$$
 This completes the proof of the result.
 \end{pf}

 \noindent
 Let $T_1$ and $T_2$ be two sets of vertices such that
  $$T_1=\Big\{v_j\in T\,|\,v_1v_j\in E(G)\Big\},\,T_2=T\backslash T_1.$$
 If $|T_1|=t_1$ and $|T_2|=t_2$, then we have
 \begin{eqnarray}
    n-k=t_1+t_2.\label{ah1}
 \end{eqnarray}
 \begin{lemma} \label{1a3} Let $G$ $(\ncong K_{1,\,n-1})$ be a connected graph of order $n$ with $m$ edges. If there is not any edge in the subgraph induced by $S$ in $G$ and $\Delta_2<\frac{2m}{n}$,
 then $\mu_2(G)>\frac{2m}{n}$.
 \end{lemma}

 \begin{pf} We have that $S$ is an independent set. Since $G$ is connected and $G\ncong K_{1,\,n-1}$, we have $\Delta_2\geq 2$. Again since $\Delta_2<\frac{2m}{n}$, by Lemma \ref{1a1},
 we have that (\ref{mo1}), (\ref{mo2}) and (\ref{mo3}) hold. From $|S_2|\geq 0$, we consider the following two cases:

 \vspace*{3mm}

 \noindent
 ${\bf Case\,1:}$ $|S_2|\geq 1$. Let $v_k$ be any one vertex in $S_2$. Then the degree of $v_k$ is $\Delta_2$. Since $S$ is an independent set, one can easily see that $K_{1,\,\Delta_2}\cup K_{1,\,\Delta_2}$
 is a subgraph of $G$ with one star $K_{1,\,\Delta_2}$ induced by $\{v_1\}$ with any $\Delta_2$ vertices in $S_1$ and the other star $K_{1,\,\Delta_2}$ induced by $\{v_k\}\cup N(v_k)$,
 where $N(v_k)\subseteq T$, $|N(v_k)|=\Delta_2$ (see, Fig. 2), that is, $K_{1,\,\Delta_2}\cup K_{1,\,\Delta_2}\subseteq G$ and hence
        $$\mu_2(G)\geq \mu_2(K_{1,\,\Delta_2}\cup K_{1,\,\Delta_2})=\Delta_2+1>\frac{2m}{n},\,\mbox{ by (\ref{mo4})}.$$

 \vspace*{3mm}

 \noindent
  ${\bf Case\,2:}$ $|S_2|=0$. From the definitions, we have  $\Delta_1=k-1+t_1$. Using this with (\ref{mo2}) and (\ref{ah1}),  we get
    \begin{equation}
  t_2<  k-1-\Delta_2. ~~~~~~\label{ah4}
   \end{equation}
 Let $G_T$ be the subgraph of $G$ induced by vertex set $T$ with  $|E(G_T)|=m_T\geq 0$. Since $S=S_1$ is an independent set, we have
 $$m=\Delta_1+(k-1)(\Delta_2-1)+m_T\leq \sum_{v_i \in T}d_i+k-1.$$
 Since $m_T\geq 0$, from the above with $\Delta_1=k-1+t_1$, we get
  \begin{eqnarray}
 (k-1)(\Delta_2-1) &\leq &  \sum_{v_i \in T}d_i-t_1 \nonumber\\[2mm]
              &= &  \sum_{v_i \in T_1}(d_i-1)+\sum_{v_i \in T_2}d_i \nonumber\\[2mm]
              &\leq & t_1(\Delta_2-2)+t_2(\Delta_2-1)    ~~~~~\mbox{ as~}~d_i\leq \Delta_2-1~\mbox{ for }~v_i\in T \nonumber\\[2mm]
              &\leq & (t_1+t_2)(\Delta_2-1). \label{km2}
 \end{eqnarray}

  \begin{center}
 \includegraphics[height=7cm,keepaspectratio]{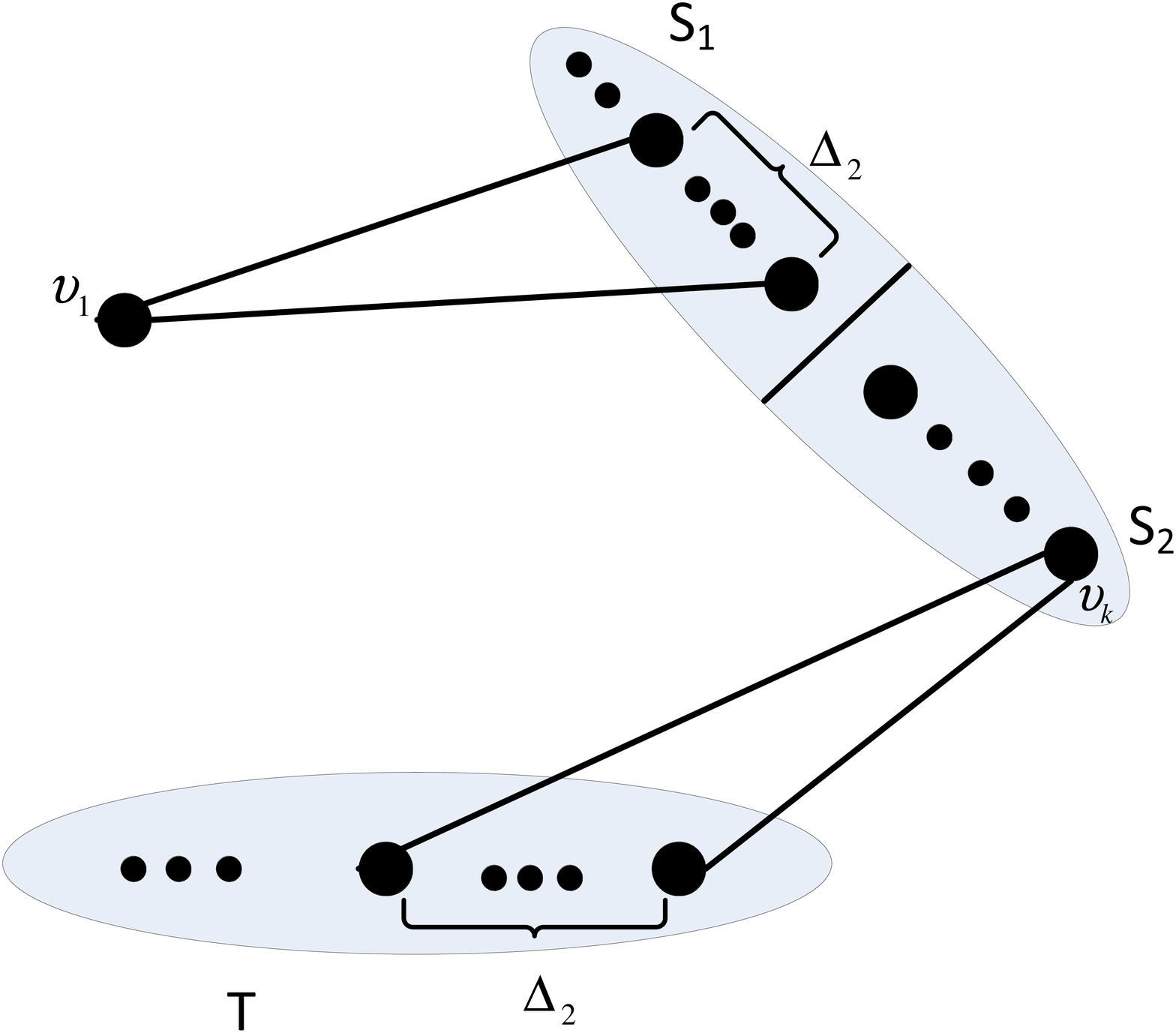}\\
  $G$\\
 \vspace*{3mm}
 \noindent
 Fig. 2. The graph $G$ contains two vertex-disjoint stars $K_{1,\,\Delta_2}$ with central vertices $v_1$ and $v_k$\,.
 \end{center}

 Since $\Delta_2\geq 2$, therefore from the above $k-1\leq t_1+t_2$. This result with (\ref{ah4}), we have $t_2<t_1+t_2-\Delta_2$, that
 is, $t_1>\Delta_2$. If $t_1(\Delta_2-2)+t_2(\Delta_2-1)=(t_1+t_2)(\Delta_2-1)$ (see, the last inequality (\ref{km2})), then $t_1=0$, a contradiction as $t_1>\Delta_2\geq 2$.
 Hence the last inequality (\ref{km2})
 is strict, that is, $k\leq t_1+t_2$ and hence $k\leq n/2$, by (\ref{ah1}). Using $k\leq t_1+t_2$ in (\ref{ah4}), we get
 $t_2< t_1+t_2-\Delta_2-1\,,~\mbox{ that is, }~t_1\geq \Delta_2+2$. Using this and (\ref{mo3}), from $\Delta_1=k-1+t_1$, we get
 $\Delta_1\geq 2\Delta_2+3.$ Also we have
 $$ \frac{2m}{n}= \frac{\Delta_1+(k-1)\Delta_2+\sum^{n}_{i=k+1}d_i}{n}\leq \Delta_2+\frac{\Delta_1-\Delta_2+k-n}{n}.$$
 From the above with $k\leq n/2$, one can easily get
    \begin{equation}\label{ah9}
  \frac{2m}{n}< \Delta_2+\frac{1}{2}.
   \end{equation}

 \vspace{3mm}

 \noindent
 We now consider the following two subcases:

 \vspace{3mm}

 \noindent
 ${\bf Subcase\,2.1:}$ There exists at least one vertex in $S$ non-adjacent to any vertex in $T_2$. Suppose that vertex $v_i$ in $S$ is not adjacent to any vertex in $T_2$.
 Then the vertex $v_i$ is adjacent to $v_1$ and $\Delta_2-1$ vertices in $T_1$ as the degree of $v_i$ is $\Delta_2$. Again since
 vertex $v_1$ is adjacent to all the vertices in $T_1$, one can easily see that $K_2\vee \overline{K}_{\Delta_2-1}$ (see, Fig. 3) is a subgraph of $G$.
 Then we have $\mu_2(G)\geq \mu_2(K_2\vee \overline{K}_{\Delta_2-1})= \Delta_2+1 > \frac{2m}{n}$, by (\ref{ah9}).

 \begin{center}
 \includegraphics[height=4cm,keepaspectratio]{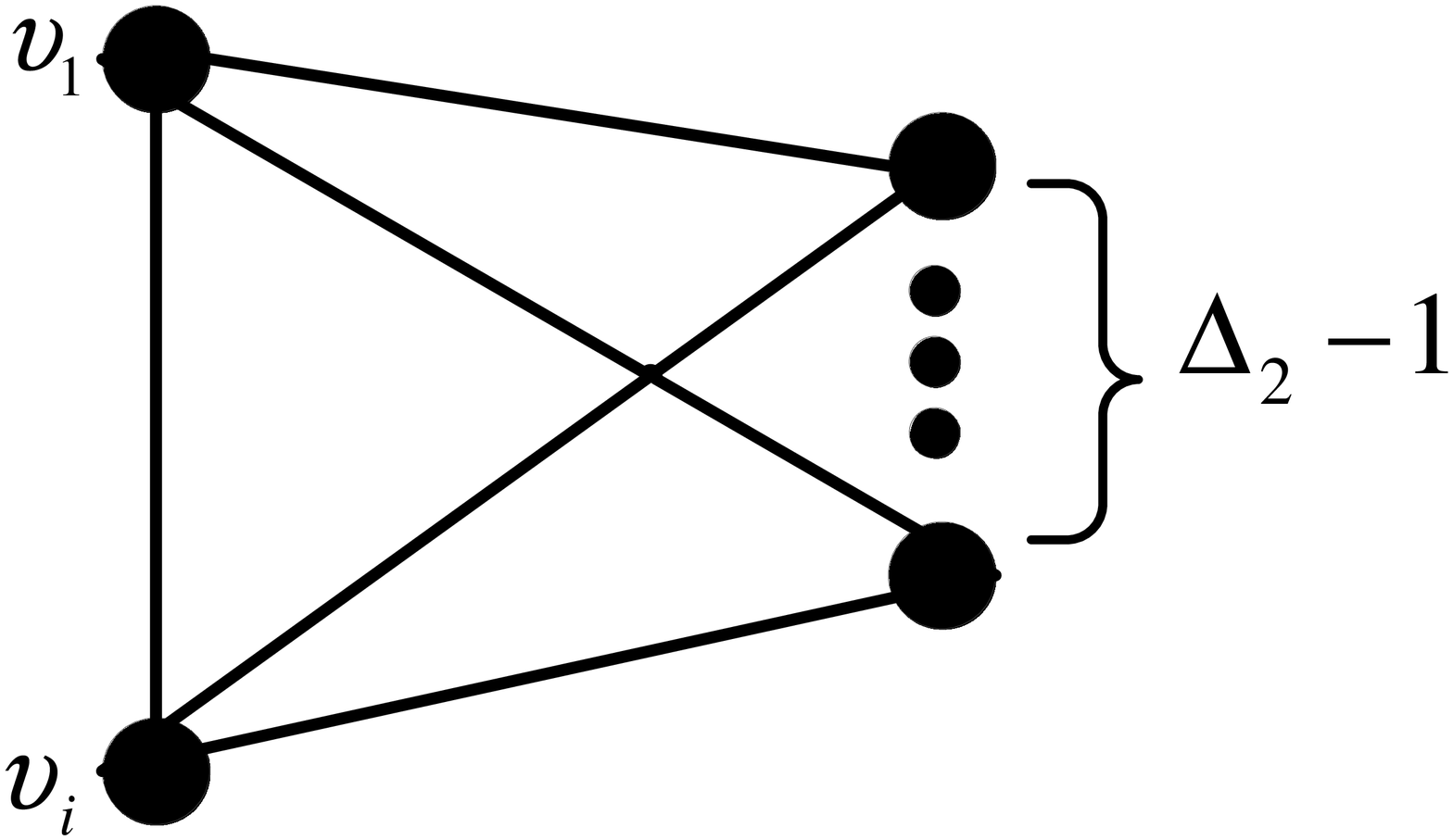}\\
 $G\cong K_2\vee \overline{K}_{\Delta_2-1}$\\
  \vspace*{3mm}
 \noindent
 Fig. 3. The graph $G\cong K_2\vee \overline{K}_{\Delta_2-1}$\,.
 \end{center}

 \vspace{3mm}

 \noindent
 ${\bf Subcase\,2.2:}$ Any one vertex in $S$ is adjacent to at least one vertex in $T_2$.\\
 We have $d_{\ell}\leq \Delta_2-1$ for any $v_{\ell}\in T_2$. Therefore
 \begin{equation}\label{ah20}
 t_2\geq \left\lceil\frac{k-1}{\Delta_2-1}\right\rceil.
 \end{equation}

 \vspace{3mm}

 If $\Delta_2=2$, then from the above, $t_2\geq k-1$. From (\ref{ah4}), we have $t_2<k-3$, a contradiction. Hence $\Delta_2\geq 3$.
 Now we have to prove that $\mu_2(G)> \frac{2m}{n}$ for $\Delta_2\geq 3$.
 Let $X=\{v_i\in T|\,d_i\leq \Delta_2-2\}$ and $|X|=x$. Then
 $$n\Delta_2< 2m\leq \Delta_1 + (k-1)\Delta_2+x(\Delta_2-2)+(n-k-x)(\Delta_2-1),$$
 that is,
 \begin{equation}
 x< \Delta_1-\Delta_2-n+k\,.\label{112}
 \end{equation}

 \vspace{3mm}

 \noindent
 Let $T'=T\backslash X$. We define two edge sets $E_{SX}$ and $E_{ST'}$, where
 $$E_{SX}=\Big\{ v_iv_j\in E(G) \,|\, v_i\in S , v_j\in X\Big\}~\mbox{ and }~E_{ST'}=\Big\{ v_iv_j\in E(G) \,|\, v_i\in S , v_j\in T'\Big\}\,.$$
 If $|X\cap T_1|=x_1$ and $|X\cap T_2|=x_2$, then
 \begin{equation}
 |E_{SX}|\leq (\Delta_2-3)x_1+(\Delta_2-2)x_2  \label{max}
 \end{equation}
 as all the vertices in $T_1$ are adjacent to $v_1$. Note that $$\,\sum_{v_i\in S}d_i=(k-1)\Delta_2=k-1+|E_{SX}|+|E_{ST'}|.$$
 From the above result with (\ref{112}) and (\ref{max}), we have
 \begin{eqnarray}
 |E_{ST'}|&=&(k-1)(\Delta_2-1)-|E_{SX}|\nonumber\\[2mm]
         &\geq&(k-1)(\Delta_2-1)-\Big((\Delta_2-3)x_1+(\Delta_2-2)x_2\Big)\nonumber\\[2mm]
         & \geq &  (k-1)(\Delta_2-1)-x_1 (\Delta_2-1)-x_2(\Delta_2-1)\nonumber\\[2mm]
         &= &(k-1-x)(\Delta_2-1) ~~~~~~~~~~~~~~\mbox{as}~~x=x_1+x_2 \nonumber\\[2mm]
         & >& \Big(n-\Delta_1+\Delta_2-1\Big)(\Delta_2-1)~~~~~~~\mbox{as}~~\Delta_2\ge 2\,.\label{mo5}
 \end{eqnarray}

 \vspace{3mm}

 \noindent
 {\bf Claim.} $\Big(n-\Delta_1+\Delta_2-1\Big)(\Delta_2-1)>k-1.$

 \vspace*{3mm}

 \noindent
 {\bf Proof of the Claim.} By contradiction, we assume that $(n-\Delta_1+\Delta_2-1)(\Delta_2-1)\leq k-1.$
  From the above and using  (\ref{ah20}), we have
  $$n-\Delta_1+\Delta_2-1\leq \frac{k-1}{\Delta_2-1}\leq t_2\,.$$
  Since $n=\Delta_1+1+t_2$, from the above, we get $\Delta_2\leq 0$, a contradiction. Hence the proof of the {\bf Claim} is finished.

 \vspace*{3mm}

 \noindent
 From the {\bf Claim} with (\ref{mo5}), we get $|E_{ST'}|/(k-1) >1.$
 Suppose that all the vertices in $S$ are adjacent to at most  one vertex in $T^{\prime}$. Then $k-1<|E_{ST'}|\leq k-1$, a
 contradiction. Hence we conclude that there is at least one vertex in $S$ (say, $v_2$) adjacent to at least two vertices (say, $v_i,\,v_j\in N(v_2)$, $i,\,j\neq 1$) of degrees $\Delta_2-1$ in $T^{\prime}$.
 We now consider the following three subcases:

 \vspace{3mm}

 \noindent
 ${\bf Subcase\,2.2.1:}$  $v_i,\, v_j\in T_1.$ Here we consider the following two subcases:

 \vspace{3mm}

 \noindent
 ${\bf Subcase\,2.2.1.1:}$  $v_iv_j \in E(G)$. In this case, $L_4$ in Lemma \ref{ap1} (Appendix), is $4\times 4$ leading principal submatrix
 corresponding to vertices $v_1$, $v_2$,  $v_i$ and $v_j$ of Laplacian matrix $L(G)$. By Lemmas \ref{k0} and \ref{ap1} with (\ref{ah9}), we conclude that
 $$\mu_2(G)\geq \mu_2(L_4)>\Delta_2+\frac{1}{2} > \frac{2m}{n}.$$

 \vspace{3mm}

 \noindent
 ${\bf Subcase\,2.2.1.2:} $ $v_iv_j \notin E(G)$.  In this case, $L_5$ in Lemma \ref{ap1} (Appendix) is $4\times 4$ leading principal submatrix corresponding
 to vertices $v_1$, $v_2$, $v_i$ and $v_j$ of Laplacian matrix $L(G)$. By Lemmas \ref{k0} and \ref{ap1} with (\ref{ah9}), we conclude that
 $$\mu_2(G)\geq \mu_2(L_5)\geq \Delta_2+1 > \frac{2m}{n}.$$

 \vspace{3mm}

 \noindent
 ${\bf Subcase\,2.2.2:}$  $v_i\in T_1$ and $v_j\in T_2.$ Here we consider the following two subcases:

 \vspace{3mm}

 \noindent
 ${\bf Subcase\,2.2.2.1:} $ $v_iv_j \in E(G)$. In this case, $L_6$ in Lemma \ref{ap1} (Appendix) is  $4\times 4$ leading principal submatrix
 corresponding to vertices $v_1$, $v_2$, $v_i$ and $v_j$ of Laplacian matrix $L(G)$. Similarly, as above, we have
 $$\mu_2(G)\geq \mu_2(L_6)> \Delta_2+\frac{1}{2}>\frac{2m}{n}.$$

 \vspace{3mm}

 \noindent
 ${\bf Subcase\,2.2.2.2:} $ $v_iv_j \notin E(G)$. In this case, $L_7$ in Lemma \ref{ap1} (Appendix) is  $4\times 4$ leading principal submatrix
 corresponding to vertices $v_1$, $v_2$, $v_i$ and $v_j$ of Laplacian matrix $L(G)$. Similarly, as above, we have
  $$\mu_2(G)\geq \mu_2(L_7)> \Delta_2+\frac{1}{2} > \frac{2m}{n}.$$

 \vspace{6mm}

 \noindent
 ${\bf Subcase\,2.2.3:}$ $v_i\,,\,v_j\in T_2.$ Here we consider the following two subcases:

 \vspace{3mm}

 \noindent
 ${\bf Subcase\,2.2.3.1:}$ $v_iv_j \in E(G)$: In this case, $L_8$ in Lemma \ref{ap1} (Appendix) is  $4\times 4$ leading principal submatrix
 corresponding to vertices $v_1$, $v_2$, $v_i$ and $v_j$ of Laplacian matrix $L(G)$. Similarly, as above, we have
 $$\mu_2(G)\geq \mu_2(L_8)> \Delta_2+\frac{1}{2} > \frac{2m}{n}.$$

 \vspace{3mm}

 \noindent
 ${\bf Subcase\,2.2.3.2:} $ $v_iv_j \notin E(G)$: In this case, $L_9$ in Lemma \ref{ap1} (Appendix) is $4\times 4$ leading principal submatrix
 corresponding to vertices $v_1$, $v_2$, $v_i$ and $v_j$ of Laplacian matrix $L(G)$.
 Similarly, as above, we have
  $$\mu_2(G)\geq \mu_2(L_9)> \Delta_2+\frac{1}{2} > \frac{2m}{n}.$$ This completes the proof of the lemma.
 \end{pf}

 \noindent
 We are now ready to give a lower bound on $\mu_2$ of connected graph $G$ and characterize the extremal graphs.
  \begin{theorem} \label{k8} Let $G$ $(\ncong K_{1,\,n-1})$ be a connected graph of order $n$ with $m$ edges. Then
  \begin{equation}
  \mu_2(G)\geq \frac{2m}{n}   \label{th1}
  \end{equation}
  with equality holding if and only if $G\cong K_{n/2,n/2}$ $(n\mbox{ is even})$.
  \end{theorem}

 \begin{pf} Since $G\ncong K_{1,\,n-1}$, we have $\Delta_2\geq 2$. If $\Delta_2\geq \frac{2m}{n}$, then by (\ref{1ekine1}), we have
 \begin{equation}
 \mu_2(G)\geq \Delta_2\geq \frac{2m}{n}  \nonumber
 \end{equation}
 and (\ref{th1}) holds. Otherwise, $\Delta_2<\frac{2m}{n}$. By (\ref{mo3}), we conclude that there are at least $\Delta_2+1$ vertices of degree $\Delta_2$.
 From $(\ref{mo2})$, we get $\Delta_1\geq \Delta_2+1+|T|$. Since $S\cup T\cup \{v_1\}=V(G)$, therefore there are at least $\Delta_2+1$ vertices in $S$, adjacent to
 $v_1$. Hence $|S_1|\geq \Delta_2+1\geq 3$ and $|S_2|\geq 0$. If there is an edge in the subgraph induced by $S$ in $G$, then by Lemma
 \ref{1a2}, we get the required result in (\ref{th1}). Otherwise, there is not any edge in the subgraph induced by $S$ in
 $G$. By Lemma \ref{1a3}, again we get the required result in (\ref{th1}). The first part of the theorem is done.

 \vspace*{3mm}

 Suppose that equality holds in (\ref{th1}).
 By Lemmas \ref{1a2} and \ref{1a3}, we must have $\mu_2(G)=\Delta_2=\frac{2m}{n}$. By Lemma \ref{k2}, we have $\Delta_1=\Delta_2$ as $G\ncong K_{1,\,n-1}$.
 Thus we have $\Delta_1=\frac{2m}{n}$ and hence $G$ is a regular graph.
 Again from $n\mu_2(G)=2m$, we have
 $$\mu_1(G)-2\mu_2(G)=\sum_{i=3}^{n-1} \left(\mu_2(G)-\mu_i(G)\right)\geq 0,~\mbox{ i.e., }~\mu_1(G)\geq 2\mu_2(G)=\frac{4m}{n}.$$
 On the other hand, by Lemma \ref{kk}, we have $$\mu_1(G)\leq \Delta_1+\Delta_2=2\Delta_2=\frac{4m}{n}.$$
 From the above results, we conclude that
 $$\mu_1(G)=\frac{4m}{n}~\mbox{ and }~\mu_2(G)=\mu_3(G)=\cdots=\mu_{n-1}(G)=\frac{2m}{n}.$$
 By Lemma \ref{t4}, we get $G\cong K_{\frac{n}{2},\frac{n}{2}}$ $(n$ is even$)$ as $G\ncong K_{1,\,n-1}$ and $\mu_1(G)\neq \mu_2(G)$.

 \vspace*{3mm}

 \noindent
 Conversely, one can easily see that the equality holds in (\ref{th1}) for $K_{n/2,n/2}$ $(n\mbox{ is even})$.
 \end{pf}

 \noindent
 We now generalize our result in Theorem \ref{k8} as follows:
 \begin{theorem} \label{1mk1} Let $G\,(\ncong K_{1,\,n-r-1}\cup r\,K_1,\,0\leq r\leq \lceil\frac{n}{2}\rceil-2$ or $r=n-2)$ be a graph of order $n>2$ with $m$ edges. Then
 \begin{equation}
 \mu_2(G)\geq \frac{2m}{n}\label{ek1}
 \end{equation}
 with equality holding if and only if $G\cong n\,K_1$ or $G\cong K_{n/2,n/2}$ $(n\mbox{ is even})$ or $G\cong K_{1,\,n/2}\cup (n/2-1)\,K_1$ $(n$ is even$)$.
 \end{theorem}

 \begin{pf} If $G\,(\ncong K_{1,\,n-1})$ is a connected graph, then by Theorem \ref{k8}, we get the required result in (\ref{ek1}).
 Otherwise, $G$ is a disconnected graph. For $G\cong n\,K_1$, (\ref{ek1}) holds. Then there exists an edge in $G$. We can assume that $G_i$ is the $i$-th connected component of order $n_i$
 with $m_i>0$ edges in $G$ $(i=1,\,2,\ldots,\,k)$ such that
 $G\cong \bigcup^k_{i=1}\,G_i\cup s\,K_1$ $(s\geq 0)$. First we assume that $k=1$ and then $s>0$. If $G_1\cong K_{1,\,n-r-1}$,
 $\Big\lceil\frac{n}{2}\Big\rceil-1\leq r\leq n-3$, then
 \begin{eqnarray}
 \mu_2(G)=\mu_2(G_1)=1\geq \frac{2m}{n}~~\mbox{ as }2\leq m=n-r-1\leq \Big\lfloor\frac{n}{2}\Big\rfloor.\label{kine1}
 \end{eqnarray}
 Otherwise, $G_1\ncong  K_{1,\,n-r-1}~(1\leq r\leq n-2)$. By Theorem \ref{k8}, we have
 \begin{eqnarray}
  \mu_2(G)=\mu_2(G_1)\geq \frac{2m_1}{n_1}=\frac{2m}{n_1}>\frac{2m}{n}~\mbox{ as }n>n_1.\nonumber
  \end{eqnarray}

 \noindent
 Next we assume that $k\geq 2$. Without loss of generality we can assume that
 \begin{equation}
 \frac{2m_1}{n_1}\geq \frac{2m_2}{n_2}\geq \cdots\geq \frac{2m_k}{n_k}.\label{ac1}
 \end{equation}
 We now prove $\frac{2m}{n}\leq \frac{2m_1}{n_1}$ by contradiction. For this, we suppose that $\frac{2m}{n}>\frac{2m_1}{n_1}$. Then we have
         $$\frac{2m}{n}>\frac{2m_1}{n_1}\geq \frac{2m_2}{n_2}\geq \cdots\geq \frac{2m_k}{n_k},~\mbox{ that is, }~2m\,n_i>2n\,m_i,~i=1,\,2,\ldots,\,k,$$
     $$\mbox{ that is, }~2m\,\sum\limits^k_{i=1}\,n_i>2n\,\sum\limits^k_{i=1}\,m_i,~\mbox{ that is, }~2m(n-s)>2mn,~\mbox{ a contradiction.}$$

 \noindent
 If $G_1\ncong K_{1,\,n_1-1}$ for some $n_1\geq 2$, then by Theorem \ref{k8} we have
 \begin{equation}
 \mu_2(G)\geq \mu_2(G_1)\geq \frac{2m_1}{n_1}\geq \frac{2m}{n}.\label{ek2}
 \end{equation}
 Otherwise, $G_1\cong K_{1,\,n_1-1}$ for some $n_1\geq 2$ and we have $\mu_1(G_1)=n_1\ge 2$.
 Since $G_2$ is a connected graph of order at least 2, we have $K_2\subseteq G_2$ ($K_2$ is a subgraph of $G_2$) and therefore $\mu_1(G_2)\ge \mu_1(K_2)=2$.
 From the above, we get $$\mu_2(G)\ge 2 >\frac{2(n_1-1)}{n_1}=\frac{2m_1}{n_1}\ge \frac{2m}{n}.$$
 The first part of the theorem is proved.

 \vspace*{3mm}

 \noindent
 Suppose that equality holds in (\ref{ek1}). Then all inequalities in the above argument must be equalities.
 For connected graph $G$, by Theorem \ref{k8}, $G\cong K_{n/2,n/2}$ $(n\mbox{ is even})$. For disconnected graph $G$,
 the equality holds in (\ref{kine1}) and (\ref{ek2}). From the equality in (\ref{kine1}), we get
    $$\mu_2(G)=\mu_2(G_1)=1=\frac{2m}{n},~~2\leq m=n-r-1\leq \Big\lfloor\frac{n}{2}\Big\rfloor,$$
 that is, $m=n-r-1=n/2$. In this case $G\cong K_{1,\,n/2}\cup (n/2-1)\,K_1$ ($n$ is even).

 \noindent
 From the equality in (\ref{ek2}), we get
 \begin{equation}
 \mu_2(G)=\mu_2(G_1)=\frac{2m_1}{n_1}=\frac{2m}{n}.\label{ad1}
 \end{equation}
 We now consider the following two cases:

 \noindent
 ${\bf Case\,1:}$ $\frac{2m_1}{n_1}=\frac{2m_k}{n_k}$. Then
 $$\mu_2(G)=\mu_2(G_1)=\frac{2m_1}{n_1}=\frac{2m_2}{n_2}=\cdots=\frac{2m_k}{n_k}=\frac{2m}{n}.$$
 From the second equality, by Theorem \ref{k8}, we have $G_1\cong
 K_{n_1/2,\,n_1/2}$. Since $\mu_1(G_1)=n_1>\frac{2m_1}{n_1}=\frac{2m}{n}=\mu_2(G)$, from the above, we must have
   $$\mu_1(G_i)\leq \mu_2(G)=\frac{2m_i}{n_i},~i=2,\,3,\ldots,\,k.$$
 But we have
    $$\mu_1(G_i)\geq \Delta(G_i)+1>\frac{2m_i}{n_i},~i=2,\,3,\ldots,\,k,$$
 which gives a contradiction.

 \vspace*{3mm}

 \noindent
 ${\bf Case\,2:}$ $\frac{2m_1}{n_1}>\frac{2m_k}{n_k}$. From (\ref{ac1}), we have
       $$2m_i\leq \frac{2m_1}{n_1}\,n_i,~~i=2,\,3,\ldots,\,k-1.$$
 Thus we have
    $$2m=2\,\sum\limits^k_{i=1}\,m_i<\frac{2m_1}{n_1}\,\sum\limits^k_{i=1}\,n_i=\frac{2m_1}{n_1}\,(n-s),~\mbox{i.e., }~\frac{2m}{n}\leq \frac{2m}{n-s}<\frac{2m_1}{n_1},$$
 a contradiction by (\ref{ad1}).

 \vspace*{2mm}

 \noindent
 Conversely, one can easily see that the equality holds in (\ref{ek1}) for $n\,K_1$ or $K_{n/2,n/2}$ $(n\mbox{ is even})$ or $K_{1,\,n/2}\cup (n/2-1)\,K_1$.
 \end{pf}

 \vspace*{3mm}

 \noindent
 We now give an upper bound on the third smallest Laplacian eigenvalue $\mu_{n-2}(G)$ in terms of $m$ and $n$.
\begin{theorem} Let $G\,(\ncong (K_1\cup K_{n-r-1})\vee K_r,\,0\leq r\leq \lceil\frac{n}{2}\rceil-2$ or $r=n-2)$. Then
 \begin{equation}
 \mu_{n-2}(G)\leq \frac{2m}{n}+1  \label{t2}
 \end{equation}
 with equality holding if and only if $G\cong K_n$ or $G\cong 2\,K_{n/2}$ $(n\mbox{ is even})$ or $G\cong (K_1\cup K_{\frac{n}{2}})\vee K_{\frac{n}{2}-1}$ $(n\mbox{ is even})$.
 \end{theorem}

 \begin{pf} Since $G\ncong (K_1\cup K_{n-r-1})\vee K_r,\,0\leq r\leq \lceil\frac{n}{2}\rceil-2,\,r=n-2$, we have $\overline{G}\ncong K_{1,\,n-r-1}\cup r\,K_1$,
 $0\leq r\leq \lceil\frac{n}{2}\rceil-2,\,r=n-2$. Then, by Theorem \ref{1mk1} and Lemma \ref{j2}, we have
   $$\mu_2(\overline{G})\geq \frac{2\overline{m}}{n},~\mbox{ that is, }~n-\mu_{n-2}(G)\geq \frac{n\,(n-1)-2m}{n},~\mbox{ that is, }~\mu_{n-2}(G)\leq \frac{2m}{n}+1.$$

 \vspace*{3mm}

 \noindent
 By Theorem \ref{1mk1}, the equality holds in (\ref{t2}) if and only if  $\overline{G}\cong n\,K_1$ or $\overline{G}\cong K_{n/2,n/2}$ $(n\mbox{ is even})$
 or $\overline{G}\cong K_{1,\,n/2}\cup (n/2-1)\,K_1$ ($n$ is even), that is, $G\cong K_n$ or $G\cong 2\,K_{n/2}$ $(n\mbox{ is even})$ or $G\cong (K_1\cup K_{\frac{n}{2}})\vee K_{\frac{n}{2}-1}$
 ($n$ is even).
 \end{pf}

 \section{Solution of Problem \ref{prob3} and application to Laplacian energy}

 In this section, we provide the answer to Problem \ref{prob3}. Using this solution, we give an upper bound for Laplacian energy of graphs.

 \vspace*{3mm}

 \noindent
 By the definition of $\sigma(G)$, we can rewrite Theorem \ref{1mk1} as follows:
 \begin{theorem} \label{k11} Let $G\,(\ncong K_{1,\,n-r-1}\cup r\,K_1,\,0\leq r\leq \lceil\frac{n}{2}\rceil-2$ or $r=n-2)$ be a graph of order $n>2$. Then $\sigma(G)\geq 2.$
 \end{theorem}

 \noindent
 If $G\cong K_{1,\,n-r-1}\cup r\,K_1\,(0\leq r\leq \lceil\frac{n}{2}\rceil-2$ or $r=n-2)$, then one can easily see that $\mu_2(G)<\frac{2m}{n}$.
 This result with Theorem \ref{k11} leads to the following result which is a complete solution to Problem \ref{prob3}.
 \begin{theorem}\label{sigma1} Let $G$ be a graph of order $n$. Then $\sigma=1$ if and only if $G\cong K_{1,\,n-r-1}\cup r\,K_1\,(0\leq r\leq \lceil\frac{n}{2}\rceil-2$ or $r=n-2)$.
 \end{theorem}

 \vspace*{3mm}

 \noindent
 Using the above result, we can improve an upper bound on Laplacian energy of graphs. The following upper bound on Laplacian energy of graphs is obtained in \cite{KI-MO}:
 \begin{theorem} \label{2k6} Let $G$ be a graph of order $n$ with $m$ $(\geq \frac{n}{2})$ edges and maximum degree $\Delta_1$. Then
 \begin{equation}
  LE(G)\leq 4m-2\Delta_1-\frac{4m}{n}+2.\label{11e0}
 \end{equation}
 \end{theorem}

 \noindent
 Now we improve this upper bound in the following:
 \begin{theorem} \label{k13} Let $G\,(\ncong K_{1,\,n-r-1}\cup r\,K_1,\,0\leq r\leq \lceil\frac{n}{2}\rceil-2$ or $r=n-2)$ be a graph of order $n$ with $m$ $(\geq \frac{n}{2})$ edges
 and maximum degree $\Delta_1$. Then
 \begin{equation}
  LE(G)\leq 4m-2\Delta_1-\frac{8m}{n}+4.\label{new}
 \end{equation}
  \end{theorem}

 \begin{pf} By Theorem \ref{k11}, we have $\sigma\geq 2$. From Lemma \ref{t1}, one can easily see that
 $$\sum^\sigma_{i=1}\lambda_i(L(G))\leq  \sum^\sigma_{i=1}\lambda_i\left(L(K_{1,\,\Delta_1})\right)+\sum^\sigma_{i=1}\lambda_i\left(L(G\backslash K_{1,\,\Delta_1}\right)),$$
 where $\Delta_1$ is the maximum degree of $G$ (For subgraph $H$ of $G$, let $G\backslash H$ be a subgraph of $G$ such that $V(G\backslash H)=V(G)$ and $E(G\backslash H)=E(G)\backslash E(H)$). Since
 \begin{eqnarray}
 \sum^\sigma_{i=1}\lambda_i\left(L(K_{1,\,\Delta_1})\right)=\sum^\sigma_{i=1}\,\mu_i(K_{1,\,\Delta_1})\leq \Delta_1+\sigma\,\mbox{ and }\,\sum^\sigma_{i=1}\lambda_i\left(L(G\backslash
 K_{1,\,\Delta_1})\right)&=&\sum^\sigma_{i=1}\,\mu_i(G\backslash K_{1,\,\Delta_1})\nonumber\\[2mm]
 &\leq& 2(m-\Delta_1),\nonumber
 \end{eqnarray}

 \noindent
 from the above, we get
 \begin{equation}
 S_{\sigma}(G)=\sum^\sigma_{i=1}\mu_i(G)\leq 2m-\Delta_1+\sigma.\nonumber
 \end{equation}

 \vspace*{3mm}

 \noindent
 Using the above result in (\ref{ghs1}), we get
 \begin{eqnarray}
 LE(G)=2S_{\sigma}(G)-\frac{4m\sigma}{n}\leq 4m-2\Delta_1-2\sigma\left(\frac{2m}{n}-1\right),\nonumber
 \end{eqnarray}
 which gives the required result in (\ref{new}) by $\sigma\geq 2$ and $m\geq \frac{n}{2}$.
 \end{pf}

 \begin{remark} \label{k13} For  graph  $G\ncong K_{1,\,n-r-1}\cup r\,K_1\,(0\leq r\leq \lceil\frac{n}{2}\rceil-2$ or $r=n-2)$ with $m\geq \frac{n}{2}$, the upper bound in $(\ref{new})$
 is always better than the upper bound in $(\ref{11e0})$.
 \end{remark}

 \vspace*{4mm}

 \noindent
 {\it Acknowledgement.} The first author was supported by the Sungkyun research fund, Sungkyunkwan University, 2017, and National Research Foundation of the Korean government with grant No.
 2017R1D1A1B03028642, and the second author was supported by the National Research Foundation of Korea (NRF) grant funded by the Korea government (MSIP)
 (2016R1A5A1008055).

 \section{Appendix}

 \begin{lemma} \label{ap1}
 Let $L_i$ $(4\leq i\leq 9)$ be the square matrices of order $4$ as follows:
  \[L_4=\left[ \begin{array}{cccc}
 \Delta_1 & -1 & -1& -1 \\
 -1 &\Delta_2 & -1& -1 \\
  -1 &-1 & \Delta_2-1& -1 \\
  -1 &-1 & -1& \Delta_2-1
 \end{array} \right],~~ L_5= \left[ \begin{array}{cccc}
 \Delta_1 & -1 & -1& -1 \\
 -1 &\Delta_2 & -1& -1 \\
  -1 &-1 & \Delta_2-1& 0 \\
  -1 &-1 & 0& \Delta_2-1
 \end{array} \right],\]\vspace*{2mm}
\[L_6= \left[ \begin{array}{cccc}
 \Delta_1 & -1 & -1& 0 \\
 -1 &\Delta_2 & -1& -1 \\
  -1 &-1 & \Delta_2-1& -1 \\
  0 &-1 &-1& \Delta_2-1
 \end{array} \right],~~L_7= \left[ \begin{array}{cccc}
 \Delta_1 & -1 & -1& 0 \\
 -1 &\Delta_2 & -1& -1 \\
  -1 &-1 & \Delta_2-1& 0 \\
  0 &-1 & 0& \Delta_2-1
 \end{array} \right],\]\vspace*{2mm}
 \[L_8= \left[ \begin{array}{cccc}
 \Delta_1 & -1 & 0& 0 \\
 -1 &\Delta_2 & -1& -1 \\
  0 &-1 & \Delta_2-1&-1 \\
  0 &-1 & -1& \Delta_2-1
 \end{array} \right],~~L_9=\left[ \begin{array}{cccc}
 \Delta_1 & -1 & 0& 0 \\
 -1 &\Delta_2 & -1& -1 \\
  0 &-1 & \Delta_2-1&0 \\
  0 &-1 & 0& \Delta_2-1
 \end{array} \right],\]

 \vspace*{2mm}
 \noindent
 where $\Delta_1\ge 2\Delta_2+3$ $(\Delta_2\geq 2)$. Then
   $$\mu_2(L_i)>\Delta_2+\frac{1}{2}~\mbox{ for }~i=4,\,6,\,7,\,8,\,9,\mbox{ and }~\mu_2(L_5)\geq \Delta_2+1.$$
 \end{lemma}

\begin{pf} The characteristic polynomial of $L_4$ is
  \begin{eqnarray}
 f(x)=x^4& +& \left(-\Delta_1 - 3\Delta_2 + 2\right)x^3 + \left(3\Delta_1 \Delta_2 + 3\Delta_2^2 - 2\Delta_1 - 4\Delta_2 - 5\right)x^2
 \nonumber\\[2mm]
        &+& \Big(-3\Delta_1 \Delta_2^2 - \Delta_2^3 + 4\Delta_1 \Delta_2 + 2\Delta_2^2 + 2\Delta_1 + 8\Delta_2 + 2\Big)x \nonumber\\[2mm]
        &+& \Delta_1 \Delta_2^3 -2\Delta_1 \Delta_2^2 - 2\Delta_1 \Delta_2 - 3\Delta_2^2 - 2\Delta_2.\nonumber
 \end{eqnarray}
 Using Lemma \ref{k0} with $\Delta_1\ge 2\Delta_2+3$, we have $\mu_1(L_4)\geq \Delta_1 \geq 2\Delta_2+3.$ On the other hand we have
 $$f\left(\Delta_2+1\right)=-\Delta_1+\Delta_2<0 ~~~  \mbox{and }~~~
 f\left(\Delta_2+\frac{1}{2}\right)=\frac{3}{8}\Delta_1-\frac{3}{8}\Delta_2+\frac{1}{16}>0\,.$$
 Hence  $\mu_2(L_4)>\Delta_2+\frac{1}{2}$.

 \vspace{3mm}

 \noindent
 The characteristic polynomial of $L_5$ is
 \begin{eqnarray}
 f(x)=x^4 &+& (-\Delta_1 - 3\Delta_2 + 2)x^3 + (3\Delta_1 \Delta_2 + 3\Delta_2^2 - 2\Delta_1 - 4\Delta_2 - 4)x^2 \nonumber\\[2mm]
          &+&\Big(-3\Delta_1 \Delta_2^2 - \Delta_2^3 + 4\Delta_1 \Delta_2 + 2\Delta_2^2 + \Delta_1 + 7\Delta_2 - 2\Big)x \nonumber\\[2mm]
          &+& \Delta_1 \Delta_2^3 -2\Delta_1 \Delta_2^2 - \Delta_1 \Delta_2 - 3\Delta_2^2 + 2\Delta_1 + 3. \nonumber
 \end{eqnarray}
 Using Lemma \ref{k0} with  $\Delta_1\ge 2\Delta_2+3$, we have $\mu_1(L_5)\geq \Delta_1 \geq 2\Delta_2+3$. On the other hand, we have $f(\Delta_2+1)=0$.
 Hence $\mu_2(L_5)\geq \Delta_2+1.$

 \vspace{3mm}

  \noindent
 The characteristic polynomial of $L_6$ is
 \begin{eqnarray}
 f(x)=x^4 &+& (-\Delta_1 - 3\Delta_2 + 2)x^3 + (3\Delta_1 \Delta_2 + 3\Delta_2^2 - 2\Delta_1 - 4\Delta_2 - 4)x^2 \nonumber\\[2mm]
          &+&\Big(-3\Delta_1 \Delta_2^2 - \Delta_2^3 + 4\Delta_1 \Delta_2 + 2\Delta_2^2 + 2\Delta_1 + 6\Delta_2 - 1\Big)x \nonumber\\[2mm]
          &+& \Delta_1 \Delta_2^3 -2\Delta_1\Delta_2^2 - 2\Delta_1 \Delta_2 - 2\Delta_2^2 + \Delta_2 + 1. \nonumber
  \end{eqnarray}
 Using Lemma \ref{k0} with $\Delta_1\ge 2\Delta_2+3$, we have
 $$f\big(\Delta_2+1\big)=-\Delta_1+\Delta_2-1<0  ~~~\mbox{and}~~~
 f\left(\Delta_2+\frac{1}{2}\right)=\frac{3}{8}\Delta_1-\frac{3}{8}\Delta_2-\frac{3}{16}>0.$$  Hence
 $\mu_2(L_6)>\Delta_2+\frac{1}{2}$.

\vspace{3mm}

 \noindent
 The characteristic polynomial of $L_7$ is
  \begin{eqnarray}
 f(x)=x^4 &+& (-\Delta_1 - 3\Delta_2 + 2)x^3 + (3\Delta_1 \Delta_2 + 3\Delta_2^2 - 2\Delta_1 - 4\Delta_2 - 3)x^2\nonumber\\[2mm]
          &+& \Big(-3\Delta_1 \Delta_2^2 - \Delta_2^3 + 4\Delta_1 \Delta_2 + 2\Delta_2^2 + \Delta_1 + 5\Delta_2 - 3\Big)x \nonumber\\[2mm]
          &+& \Delta_1 \Delta_2^3 -2\Delta_1 \Delta_2^2 - \Delta_1 \Delta_2 - 2\Delta_2^2 + 2\Delta_1 + \Delta_2 + 2. \nonumber
 \end{eqnarray}
 Using Lemma \ref{k0} with $\Delta_1\ge 2\Delta_2+3$, we have $\mu_1(L_7)\geq \Delta_1 \geq 2\Delta_2+3.$ On the other hand, we have
 $$f\big(\Delta_2+1\big)=-1<0  ~~~\mbox{and} ~~~  f\left(\Delta_2+\frac{1}{2}\right)=\frac{15}{8}\Delta_1-\frac{15}{8}\Delta_2+\frac{1}{16}>0.$$
 Hence $\mu_2(L_7)>\Delta_2+\frac{1}{2}$.

 \vspace{3mm}

  \noindent
 The characteristic polynomial of $L_8$ is
   \begin{eqnarray}
  f(x)=x^4 &+& (-\Delta_1 - 3\Delta_2 + 2)x^3 + (3\Delta_1\Delta_2 + 3\Delta_2^2 - 2\Delta_1 - 4\Delta_2 - 3)x^2\nonumber\\[2mm]
           &+& \Big(-3\Delta_1 \Delta_2^2 - \Delta_2^3 + 4\Delta_1 \Delta_2 + 2\Delta_2^2 + 2\Delta_1 + 4\Delta_2 - 2\Big)x \nonumber\\[2mm]
           &+& \Delta_1 \Delta_2^3 -2\Delta_1 \Delta_2^2 - 2\Delta_1 \Delta_2 - \Delta_2^2 + 2\Delta_2. \nonumber
  \end{eqnarray}
 Using Lemma \ref{k0} with $\Delta_1\ge 2\Delta_2+3$, we have $\mu_1(L_8)\geq \Delta_1 \geq 2\Delta_2+3.$ On the other hand, we have
 $$f(\Delta_2+1)=-\Delta_1+\Delta_2-2<0 ~~~\mbox{and}~~~  f\left(\Delta_2+\frac{1}{2}\right)=\frac{3}{8}\Delta_1-\frac{3}{8}\Delta_2-\frac{23}{16}>0.$$
 Hence $\mu_2(L_8)>\Delta_2+\frac{1}{2}$.

\vspace{3mm}

  \noindent
 The characteristic polynomial of $L_9$ is
 \begin{eqnarray}
 f(x)=x^4 &+& (-\Delta_1 - 3\Delta_2 + 2)x^3 + (3\Delta_1 \Delta_2 + 3\Delta_2^2 - 2\Delta_1 - 4\Delta_2 - 2)x^2 \nonumber\\[2mm]
         &+&\Big(-3\Delta_1 \Delta_2^2 - \Delta_2^3 + 4\Delta_1 \Delta_2 + 2\Delta_2^2 + \Delta_1 + 3\Delta_2 - 4\Big)x \nonumber\\[2mm]
         &+& \Delta_1 \Delta_2^3 -2\Delta_1 \Delta_2^2 - \Delta_1 \Delta_2 - \Delta_2^2 + 2\Delta_1 + 2\Delta_2 - 1. \nonumber
  \end{eqnarray}
 Using Lemma \ref{k0} with $\Delta_1\ge 2\Delta_2+3$, we have $\mu_1(L_9)\geq \Delta_1 \geq 2\Delta_2+3.$ On the other hand, we have
 $$f\Big(\Delta_2+1\Big)=-4<0   ~~~\mbox{ and}~~~  f\big(\Delta_2+\frac{1}{2}\big)=\frac{15}{8}\Delta_1-\frac{15}{8}\Delta_2-\frac{51}{16}>0.$$
 Hence $\mu_2(L_9)>\Delta_2+\frac{1}{2}$.
 \end{pf}


\begin{thebibliography}{99}

 \bibitem{ANDER} W. N. Anderson, T. D. Morley, Eigenvalues of the Laplacian of a graph, Linear Multilinear Algebra 18 (1985) 141--145.

 \bibitem{DAS} K. C. Das, The largest two Laplacian eigenvalues of a graph, Linear Multilinear Algebra 52 (2004) 441--460.

 \bibitem{DAS3}  K. C. Das, A sharp upper bound for the number of spanning trees of a graph, Graphs and Combinatorics  23   (2007) 625--632.

 \bibitem{DA2} K. C. Das, I. Gutman, A. Sinan \c{C}evik, B. Zhou, On Laplacian energy, MATCH Commun. Math. Comput. Chem. 70 (2) (2013) 689--696.

 \bibitem{KI-MO} K. C. Das, S. A. Mojallal, On Laplacian energy of graphs, Discrete Math. 325 (2014) 52--64.

 \bibitem{KMG}  K. C. Das, S. A. Mojallal,  I. Gutman, On Laplacian energy in terms of graph invariants, Appl. Math. Comput. 268 (2015) 83--92.

 \bibitem{KMT} K. C. Das, S. A. Mojallal,  V. Trevisan, Distribution of Laplacian eigenvalues of graphs, Linear Algebra Appl. 508 (2016) 48--61.

 \bibitem{FA} K. Fan, On a theorem of Weyl concerning eigenvalues of linear transformations I, Proc. Nat. Acad. Sci. USA 35 (1949) 652--655.

 \bibitem{GMS} R. Grone, R. Meris, V. S. Sunder, The Laplacian spectrum of a graph, SIAM J. Matrix Anal. Appl. 11 (1990) 218--238.

 \bibitem{MW} J.-M. Guo, S. W. Tan, A relation between the matching number and Laplacian spectrum of a graph, Linear Algebra Appl. 325 (2001) 71--74.

 \bibitem{GuZh} I. Gutman, B. Zhou, Laplacian energy of a graph, Linear Algebra Appl. 414 (2006) 29--37.

 \bibitem{HJT} S. T. Hedetniemi, D. P. Jacobs, V.  Trevisan, Domination number and Laplacian eigenvalue distribution, European Journal of Combinatorics 53 (2016)
 66--71.

 \bibitem{HE} J. V. D. Heuvel, Hamilton cycles and eigenvalues of graphs, Linear Algebra Appl. 226-228 (1995) 723--730.

 \bibitem{LP} J. S. Li, Y. L. Pan, A note on the second largest eigenvalue of the Laplacian matrix of a graph, Linear Multilinear Algebra 48 (2000) 117--121.

 \bibitem{Me2} R. Merris, The number of eigenvalues greater than two in the Laplacian spectrum of a graph, Portugal. Math. 48 (1991) 345--349.

 \bibitem{MERIS}  R. Merris,  Laplacian matrices of graphs: A survey, Linear Algebra Appl. 197, 198  (1994) 143--176.

 \bibitem{PH} Y. -L. Pan, Y.-P. Hou, Two necessary conditions for $\lambda_2(G)=d_2(G)$, Linear Multilinear Algebra 51 (1) (2003) 31--38.

 \bibitem{PG}  S. Pirzada, H. A. Ganie, On the Laplacian eigenvalues of a graph and Laplacian energy, Linear Algebra Appl. 486 (2015) 454--468.

  \bibitem{L15} M. Robbiano, R. Jim\'enez, Applications of a theorem by Ky Fan in the theory of Laplacian energy of graphs, MATCH Commun. Math. Comput. Chem. 62 (2009) 537--552.

 \bibitem{SC} J. R. Schott, {\it Matrix Analysis for Statistics\/}, Wiley, New York, 1997.

 \bibitem{WYS} Y. Wua, G. Yub, J. Shu, Graphs with small second largest Laplacian eigenvalue, European Journal of Combinatorics 36 (2014)
 190--197.

 \end{thebibliography}
 \end{document}